\documentclass[letterpaper, 10 pt, conference]{ieeeconf} 
\IEEEoverridecommandlockouts                              

\usepackage{mathtools}
\usepackage{amssymb}
\usepackage{bbold}
\usepackage[yyyymmdd,hhmmss]{datetime}
\usepackage{bm,upgreek}
\usepackage{amssymb}  
\usepackage{cite}
\usepackage{color}
\mathtoolsset{showonlyrefs,showmanualtags}
\usepackage{dsfont}
\usepackage{booktabs}
\usepackage[final]{graphicx}

\DeclareMathOperator{\VEC}{vec}

\newcommand{\Exp}[1]{\mathbb{E}\big[ #1\big]}
\newcommand{\Prob}[1]{\mathbb{P} \big(#1\big)}   
\newcommand{\SUMI}[1]{\frac{1}{n} \sum_{i=1}^n \mathbb{1} \left( #1\right)}         
\newcommand{\argmin}{\operatornamewithlimits{argmin}} 

\newcommand*{\Cdot}{\raisebox{-0.25ex}{\scalebox{1.2}{$\cdot$}}}

\usepackage[standard]{ntheorem}
\newtheorem{Assumption}{Assumption}
\newtheorem{Problem}{Problem}

\newcommand{\B}[1]{\textcolor{black}{#1}}

\title{\LARGE \bf
A Certainty Equivalence Result in Team-Optimal Control of  Mean-Field Coupled Markov Chains
}

\author{Jalal Arabneydi and Amir G. Aghdam
\thanks{ This work has been supported in part by the Natural Sciences and Engineering Research Council of Canada (NSERC) under Grant RGPIN-262127-17, and in part by Concordia University under Horizon Postdoctoral Fellowship.}  
\thanks{Jalal Arabneydi and Amir G. Aghdam are with the  Department of Electrical and Computer Engineering, 
        Concordia University, 1455 de Maisonneuve Blvd, Montreal, QC, Canada. 
        {\tt\small Email:jalal.arabneydi@mail.mcgill.ca} and
        {\tt\small Email:aghdam@ece.concordia.ca}}%
}

\begin{document}
\maketitle

\vspace*{-5.2cm}{\footnotesize{Proceedings of IEEE Conference on Decision and Control, 2017.}}
\vspace*{4.45cm}

\thispagestyle{empty}
\pagestyle{empty}

\begin{abstract}
This paper studies a large number of homogeneous  Markov decision processes where the transition probabilities and costs are coupled in the empirical distribution of states (also called mean-field). The state of each process  is not known to others, which means that the information structure is fully decentralized.  The objective is to minimize  the average cost, defined as the empirical mean of individual costs, for which a sub-optimal solution is proposed. This solution does not depend on  the number of processes, yet  it converges to the optimal solution of the so-called mean-field sharing as the number of processes tends to infinity. Under some mild conditions, it is shown that the convergence rate of the proposed decentralized solution is proportional to the square root of the inverse of the number of processes. Finding this sub-optimal solution involves a non-smooth non-convex  optimization problem over an uncountable set, in general. To overcome this drawback,  a  combinatorial  optimization problem is introduced that achieves the same rate of convergence. 
\end{abstract}
\section{Introduction}
\subsection{Motivation}
Team-optimal control of Markov chains  have recently attracted much attention due to their potential applications in emerging areas such as  smart grids~\cite{meyn2015ancillary}, social networks~\cite{Foroutan2017social}, swarm robotics \cite{valentini2017achieving}, and transportation networks \cite{rouhieh2014optimizing}. These  applications normally involve many interconnected decision makers, wishing to collaborate in order to  minimize a common cost function~\cite{yuksel2009stochastic}.

When  the decision makers are modeled as controlled Markov chains and  joint state is known to all, the optimal solution is identified by the celebrated dynamic programming~\cite{kumar2015stochastic}.  The computational complexity of solving this dynamic program  is exponential in the number of  decision makers, in general. In addition,  at each time instant, the joint state  (a vector of  the same size as  the number of decision makers) must be communicated among all decision makers.   In practice, however,  each decision maker has limited computation and communication resources.  Due to such practical limitations, mean-field models  have received  much attention recently  for the scalablity  of their solution. Inspired  by statistical mechanics and classical physics, mean-field games were first introduced in the context of game theory in~ \cite{huang2003individual,HuangPeter2006,
huang2007large,Lasry2006I,Lasry2006II,Lasry2007mean,
weintraub2006oblivious}, and then were extended to various cases~\cite{Caines2013springer, LasryLions2011,Gomes2013survey,saldi2017markov,AdlJohWeiGol2009Oblivious}. In mean-field games, the solution concept  is \emph{Nash} strategy and the term \emph{mean-field} refers to the empirical distribution of infinite population of players. When the population is large, the effect of a single player on other players becomes negligible. Using this observation, an approximate Nash strategy is derived such that the approximation error converges to zero as the size of population goes to infinity. 


In the context of team theory, mean-field teams were first  introduced in   \cite{arabneydi2016new} and  the early results were presented  in \cite{JalalCDC2014,Jalal2015ICC,JalalCDC2015,Jalal2017linear}. In mean-field teams, the solution concept is a  \emph{team-optimal} strategy and the term \emph{mean-field} refers to the empirical distribution of finite population. In \cite{JalalCDC2014},  a dynamic programming decomposition is derived to obtain a globally optimal solution,  irrespective of the size of population (not necessarily large population), under mean-field sharing information structure. To implement the mean-field sharing, the   communication network  of agents must be connected.   In practice, however, having a connected network may not be practically feasible or economically viable, specially when the population is large. Therefore, a completely decentralized strategy is desirable in this type of problem. 

 In \cite{Nevroz2016},   a solution approach of mean-field games is adopted to find an  approximate  person-by-person optimal strategy for the finite-horizon case. The strategy is identified by a dynamic program,  and the approximation error is shown to go to zero at the rate $1 / \sqrt{n}$ as $n$ increases,  under  some Lipschitz conditions on the  dynamics, cost, and  the strategy.  In \cite{tembine2009mean},  the existence of an approximate person-by-person optimal strategy for the discounted cost infinite-horizon case is established.  The strategy is identified by an irregular Hamilton-Jacobi-Bellman equation  for  which  the solution is not necessarily the viscosity solution. Under some Lipschitz conditions, the approximation error of such a strategy is  shown to converge to zero in  distribution as $n$ increases.

In this paper, it is desired to find a completely decentralized strategy whose performance is sufficiently close to that obtained by  the mean-field sharing  strategy. Finding such a strategy is conceptually challenging because every agent has a different perspective (i.e., information) of the system and any such discrepancy would make it difficult to establish cooperation among agents.  In contrast to~\cite{Nevroz2016},  we use the dynamic program of mean-field teams   that  is fundamentally different from that of mean-field games; in addition, we do not impose any assumption on the strategy\footnote{Note that  verifying any assumption on the strategy is typically very difficult  because finding the strategy itself is an open problem, in general.}.   In contrast to~\cite{tembine2009mean}, the convergence here is in the sense of almost surely.  In contrast to both papers, we consider \textit{global} optimality rather than person-by-person optimality. 



The rest of this paper is organized as follows. In Section~\ref{sec:problem-formulation},  the problem is formulated in the context of controlled mean-field coupled Markov chains. Then, the main results are presented in Section~\ref{sec:main-results}, followed by concluding remarks given  in  Section~\ref{sec:conclusion}. 

\subsection{Notation}
Throughout the paper, $\mathbb{N}$, $\mathbb{R}_{\geq 0}$, and $\mathbb{R}_{>0}$ refer to natural numbers, non-negative real numbers, and positive real numbers, respectively. The finite set of integers $\{1,\ldots,k\}$ is denoted by $\mathbb{N}_k$. Moreover, $\Prob{\cdot}$ is the probability of a random variable; $\Exp{\cdot}$ represents the expectation of an event; $\mathbb{1}(\cdot)$ is the indicator function of a set; $\| \cdot \|_{\infty}$ represents the infinity norm of a vector, and  $| \cdot|$ denotes the absolute value of a real number or the cardinality of a set.  The short-hand notation $x_{1:t}$ is used  to denote vector $\VEC(x_1,\ldots,x_t)$. Given $n \in \mathbb{N}$ and a finite set $\mathcal{X}$,  the following spaces are defined.
\begin{table}[ht] 
\caption{Table of the spaces used in this paper.}
\setlength\tabcolsep{0em}
\begin{tabular}{l} 
  \toprule
\multicolumn{1}{c}{\textbf{Space of probability measures} }\\
\midrule
$\Delta(\mathcal{X})=\{(p_1,\ldots,p_{|\mathcal{X}|}) \big| p_k \in [0,1], k \in \mathbb{N}_{|\mathcal{X}|}, \sum_{k=1}^{|\mathcal{X}|} p_k=1\}$\\
\midrule
\multicolumn{1}{c}{\textbf{Space of empirical distributions (mean-field)} }\\
\midrule
$\mathcal{M}_n=\{(p_1,\ldots,p_{|\mathcal{X}|}) \big| p_k \in \{0, \frac{1}{n},\ldots,1\}, k \in \mathbb{N}_{|\mathcal{X}|}, \sum_{k=1}^{|\mathcal{X}|} p_k=1\}$\\
\midrule
\multicolumn{1}{c}{\textbf{Product space of unit intervals} }\\
\midrule
$\mathcal{I}(\mathcal{X})=\{(p_1,\ldots,p_{|\mathcal{X}|}) \big| p_k \in [0,1], k \in \mathbb{N}_{|\mathcal{X}|}\}$\\
\midrule
\multicolumn{1}{c}{\textbf{$\mathcal{I}(\mathcal{X})$ uniformly quantized by $\frac{1}{n}$} }\\
\midrule
$\mathcal{Q}_n=\{(p_1,\ldots,p_{|\mathcal{X}|}) \big| p_k \in \{0, \frac{1}{n},\ldots,1\}, k \in \mathbb{N}_{|\mathcal{X}|}\}$\\
  \bottomrule
\end{tabular}
\end{table}

\noindent The following relationships hold between  above spaces:
\begin{equation}
\mathcal{M}_n \subset \Delta(\mathcal{X}) \subset \mathcal{I}(\mathcal{X}) \quad  \text{and} \quad \mathcal{M}_n \subset \mathcal{Q}_n \subset \mathcal{I}(\mathcal{X}).
\end{equation}

\section{Problem formulation}\label{sec:problem-formulation}
Consider a dynamical system  consisting of  $n \in \mathbb{N}$ homogeneous agents (decision makers or controlled Markov chains)\footnote{For ease of reference, we  only use term \emph{agent} in the sequel.}  operating over a fixed finite horizon $T \in \mathbb{N}$. Let $x^i_t \in \mathcal{X}$ denote the state of agent $i \in \mathbb{N}_n$   at time $t=\mathbb{N}_T$ and $u^i_t \in \mathcal{U}$ represent its control action. Let also $m_t$ be the empirical distribution of states at time $t$, i.e.,
\begin{equation}\label{eq:mean-field}
m_t(x)=\frac{1}{n}\sum_{i=1}^n \mathbb{1}(x^i_t=x), \quad x \in \mathcal{X},
\end{equation}
where $m_t \in \mathcal{M}_n$. At time $t \in \mathbb{N}_T$, the state of agent $i \in \mathbb{N}_n$ evolves  as follows:
\begin{equation}\label{eq:model1}
x^i_{t+1}=f_t(x^i_t,u^i_t,w^i_t,m_t),
\end{equation}
where $w^i_t \in \mathcal{W}$ is the local noise  of agent~$i$ at time $t$.  The spaces $\mathcal{X}$, $\mathcal{U}$, and $\mathcal{W}$ are finite-valued\footnote{The main results of this paper hold for any measurable set $\mathcal{W}$ as long as the variance of the noise process is uniformly bounded.} and  it is assumed that the primitive random variables are defined on a common probability space. The initial states $\mathbf x_1:=\VEC(x^1_1,\ldots,x^n_1)$ and     noises  $\mathbf w_t:=\VEC(w^1_t,\ldots,w^n_t)$ are i.i.d.  random variables\footnote{Note that the variances of primitive random variables are finite when $\mathcal{X}$ and $\mathcal{W}$ are finite-valued spaces.}.  Also, $\{\mathbf x_1, \mathbf w_1,\ldots, \mathbf w_T\}$ are mutually independent.

The dynamics of agent $i \in \mathbb{N}_n$, given by \eqref{eq:model1}, may be equivalently expressed in the form  of \emph{controlled mean-field coupled  Markov chains} as follows:
\begin{equation}\label{eq:model2}
\Prob{x^i_{t+1}=y|x^i_t=x,u^i_t=u, m_t=m},
\end{equation} 
where the above expression corresponds to the probability of  the realizations of $w^i_t$ that take state $x \in \mathcal{X}$ to state $y \in \mathcal{X}$ under action $u \in \mathcal{U} $  when  the mean-field is $m \in \mathcal{M}_n$, i.e.,
\begin{equation}\label{eq:relation-models}
\Prob{y|x,u,m}=\sum_{w \in \mathcal{W}} \mathbb{1}\left(f_t(x,u,m,w)=y\right)\Prob{w^i_t=w}.
\end{equation}
For ease of display, we occasionally interchange \eqref{eq:model1} and~\eqref{eq:model2} in the sequel. Define $\mathbf x_t:=\VEC(x^1_t,\ldots,x^n_t)$ and $\mathbf u_t:=\VEC(u^1_t,\ldots,u^n_t)$. At time $t$, the system incurs a per-step cost given by
\begin{equation}
c_t(\mathbf x_t,\mathbf u_t)= \frac{1}{n} \sum_{i=1}^n \ell_t(x^i_t,u^i_t,m_t),
\end{equation}
where $\ell_t: \mathcal{X} \times \mathcal{U} \times \mathcal{M}_n \rightarrow \mathbb{R}_{\geq 0} $.  Denote by $I^i_t$  the information available to agent $i$ at time $t$. Then, 
\begin{equation}
u^i_t=g^i_t(I^i_t),
\end{equation}
where $g^i_t$ is called  \emph{the control law} of agent $i$ at time $t$.  The collection of control laws $\mathbf g:=\{g^1_t,\ldots,g^n_t\}_{t=1}^T$ is called a \emph{strategy}. The performance of   strategy $\mathbf g$ is measured by the following function
\begin{equation}\label{eq:total cost}
J(\mathbf g)=\mathbb{E}^{\mathbf{g}}\big[\sum_{t=1}^T c_t(\mathbf x_t,\mathbf u_t)\big].
\end{equation}
In \cite{JalalCDC2014},  it is assumed that every agent $i \in \mathbb{N}_n$ observes its local state $x^i_t$  as well as the mean-field $m_t$, i.e., $I^i_t=(x^i_t,m_{t})$. Under this so called \emph{mean-field sharing} information structure, a dynamic programming decomposition is derived to obtain a globally optimal solution when agents use \emph{homogeneous control laws}, i.e.,
\begin{equation}\label{eq:mf-sharing}
\hspace{-2.7cm}\textbf{Mean-Field Sharing:} \quad  u^i_t=g_t(x^i_t,m_t).
\end{equation}

There are various methods (including consensus algorithms \cite{bishop2014distributed,olfati2006belief}) to compute and communicate  the mean-field among agents; however,  the necessary condition for  all these methods to work  is to have a connected network.  Establishing such a connected network, specially for very large number of agents, may not be feasible, both practically and economically. For this reason,  we consider  a completely decentralized information structure  in  this paper. In particular,  every agent $i \in \mathbb{N}_n$ observes only its own local state $x^i_t$ and makes  the decision $u^i_t$ as follows:
\begin{equation}\label{eq:no-sharing}
\hspace{-1cm}\textbf{Fully  Decentralized Structure:} \quad u^i_t=g_t(x^i_t,  g_{1:t-1}).
\end{equation}
\begin{Problem}\label{problem:finite}
Let $J^\ast$ denote the optimal performance under mean-field sharing information structure \eqref{eq:mf-sharing}.  It is desired  to find a sub-optimal strategy $\mathbf g$, given by \eqref{eq:no-sharing}, under which the system performance $J(\mathbf{g})$ is guaranteed  to be within $\epsilon(n)$-neighborhood  of  $J^\ast$, i.e.,
\begin{equation} \label{eq:approximation error-definition}
| J(\mathbf g) - J^\ast \big| \leq \epsilon(n). 
\end{equation}
\end{Problem}

\section{Main Results}\label{sec:main-results}
In this section,  a completely decentralized strategy is proposed, as a \emph{sub-optimal} alternative to the mean-field sharing  solution, as noted in the previous  section. In particular,  it is shown that the optimality gap, given by  \eqref{eq:approximation error-definition}, converges to zero at the rate $ 1/ \sqrt{n}$ as $n$ increases. To this end,  the following assumption is made.
\begin{Assumption}\label{assumption: finite-lipschitz}
The transition probabilities and per-step costs are Lipschitz functions in mean-field. More precisely, there exist constants $K^1_t, K^2_t \in \mathbb{R}_{>0}$, $t \in \mathbb{N}_T$, such that for every $x,y \in \mathcal{X}$, $u \in \mathcal{U}$, $z_1,z_2 \in \mathcal{I}(\mathcal{X})$,
\begin{align}
\Big| \Prob{y|x,u,z_1} - \Prob{y|x,u,z_2} \Big| &\leq K^1_t\|z_1- z_2\|_{\infty},\\
\Big| \ell_t(x,u,z_1) - \ell_t(x,u,z_2) \Big| &\leq K^2_t\|z_1- z_2\|_{\infty}.
\end{align}
\end{Assumption}
\begin{remark}
It is to be noted that every polynomial function of mean-field is a Lipschitz function because mean-field is confined to the bounded interval $\mathcal{I}(\mathcal{X})$ \cite[Corollary 12.2]{Eriksson2004}. It is worth highlighting that, according to Weierstrass Approximation Theorem~\cite{perez2008survey},  any continuous function  can be uniformly approximated as closely as desired by polynomial functions.
\end{remark}
Let $\gamma_t: \mathcal{X} \rightarrow \mathcal{U}$  be  the local map from state space~$\mathcal{X}$ to action space $\mathcal{U}$ at time $t \in \mathbb{N}_T$, i.e., from~\eqref{eq:no-sharing} 
\begin{equation}\label{eq:NS-gamma}
\gamma_t(\Cdot):=g_t(\Cdot,g_{1:t-1}).
\end{equation}
According to ~\eqref{eq:no-sharing} and \eqref{eq:NS-gamma},
\begin{equation}
u^i_t=\gamma_t(x^i_t).
\end{equation}
 Denote by $\mathcal{G}$  the set of all  mappings $\gamma:\mathcal{X} \rightarrow \mathcal{U}$ and note that   $\mathcal{G}$ is a finite set of size $|\mathcal{G}|=|\mathcal{U}|^{|\mathcal{X}|}$. For every $\gamma \in \mathcal{G}$ and $z \in \mathcal{I}(\mathcal{X})$, define
\begin{equation}\label{eq:mean-field-f}
\hat f_t(z,\gamma)(\cdot):=\sum_{x \in \mathcal{X}}  z(x)  \Prob{\cdot|x,\gamma(x),z},
\end{equation}
and
\begin{equation}\label{eq:mean-field-cost}
\hat c_t(z,\gamma):=\sum_{x \in \mathcal{X}}  z(x)  \ell_t(x,\gamma(x),z).
\end{equation}
\begin{Lemma}\label{lemma:Lipschitz}
Let Assumption~\ref{assumption: finite-lipschitz} hold. Then, there exist constants $K^3_t,K^4_t \in \mathbb{R}_{>0}$, $t \in \mathbb{N}_T$,  such that for every $\gamma \in \mathcal{G}$ and  $z_1,z_2 \in~\mathcal{I}(\mathcal{X})$,
\begin{align}
\|\hat f_t(z_1,\gamma) -\hat f_t(z_2,\gamma) \|_{\infty} \leq K^3_t \|z_1 - z_2 \|_{\infty},\\
|\hat c_t(z_1,\gamma) -\hat c_t(z_2,\gamma) |\leq K^4_t \|z_1 - z_2 \|_{\infty}.
\end{align}
\end{Lemma}
\begin{Proof}
According to \cite[Theorem 12.1]{Eriksson2004} and \cite[Theorem 12.4]{Eriksson2004} any linear combination or product of Lipschitz functions is a Lipschitz function as well. Hence,  function $\hat f_t(z,\gamma)$ given by  \eqref{eq:mean-field-f} is Lipschitz because it is a linear combination of the product of two Lipschitz functions $z$ and $\Prob{.|x,\gamma(x),z}$.  Analogously, function $\hat c_t(z,\gamma)$ given by  \eqref{eq:mean-field-cost} is Lipschitz because it is a linear combination of the product of two Lipschitz functions $z$ and $\ell_t(x,\gamma(x),z)$. $\hfill \blacksquare$
\end{Proof}
\begin{Lemma}\label{lemma:square-root-iid}
Consider $n$ i.i.d. random variables $W^i \in \mathcal{W}$ with common probability mass function $P(W)$. Then, for every realization $w \in \mathcal{W}$, one has
\begin{equation} \label{eq:square-root-iid}
\Exp{\Big| \frac{1}{n} \sum_{i=1}^n \mathbb{1}(W^i=w) - P(W=w)\Big|} \leq \mathcal{O}(\frac{1}{\sqrt{n}}).
\end{equation}
\end{Lemma}
\begin{Proof}
Let $b$  be a random variable on $\mathbb{R}$ and $n \in \mathbb{R}_{>0}$. Then, as the first step, the following inequality is established
\begin{equation}\label{eq:inequality-useful}
\Exp{|b| }\leq \frac{\sqrt{n}}{2} \Exp{b^2} + \frac{1}{2 \sqrt{n}}.
\end{equation}
This   follows  immediately from the inequality $0 \leq (b \pm \frac{1}{\sqrt{n}})^2$, after rewriting it in the following form  $-(b^2+\frac{1}{n}) \leq \pm \frac{2}{\sqrt{n}} b \leq (b^2+\frac{1}{n})$ and noting that the expectation operator is monotone. 

In the second step, one has
\begin{align}\label{eq:square-1-n}
\mathbb{E} \Bigg[ \Big(&\frac{1}{n} \sum_{i=1}^n \mathbb{1}(W^i=w) - P(W=w)\Big)^2 \Bigg] \nonumber \\
&=\frac{1}{n^2} \Exp{ \Big( \sum_{i=1}^n \left(\mathbb{1}(W^i=w) - P(W=w)\right)\Big)^2 } \nonumber \\
& \substack{(a)\\=} \frac{1}{n^2}  \sum_{i=1}^n \Exp{\left(\mathbb{1}(W^i=w) - P(W=w)\right)^2}\substack{(b)\\ \leq}  \mathcal{O}(\frac{1}{n}),\quad 
\end{align}
where $(a)$ follows from the fact that the random variables $\left(\mathbb{1}(W^i=w) - P(W=w)\right)$ have zero-mean\footnote{By definition of expectation function, we have $ \Exp{\mathbb{1}(W^i=w)}=P(W=w), \forall i \in \mathbb{N}_n$.} for any $i \in \mathbb{N}_n$ and are mutually independent too,  which implies that the cross-terms are equal to zero; $(b)$ follows from the inequality  $\Exp{\left(\mathbb{1}(W^i=w) - P(W=w)\right)^2} \leq 1$, $\forall i \in \mathbb{N}_n$.  The proof is now complete by  virtue of  inequalities~ \eqref{eq:inequality-useful} and ~\eqref{eq:square-1-n}.~$\hfill \blacksquare$
\end{Proof}

\begin{Lemma}\label{lemma:bar f-mean-field}
Given $m_t \in \mathcal{M}_n$ and $\gamma_t \in \mathcal{G}$, there exists a function $\bar f_t$ such that  the dynamics of the mean-field is (almost surely) expressed as
\begin{equation}
m_{t+1} \substack{{a.s.}\\{=}} \bar{f}_t(m_t,\gamma_t,\mathbf w_t).
\end{equation}
\end{Lemma}
\begin{Proof}
For every $y \in \mathcal{X}$, $m \in \mathcal{M}_n$ and $\gamma \in \mathcal{G}$, define 
\begin{align}\label{eq:def-bar-f}
\bar f_t(m,\gamma,\mathbf w_t)(y): =& \sum_{x \in \mathcal{X}} \sum_{w \in \mathcal{W}}  \Bigg[\mathbb{1}\left( f_{t}(x),\gamma(x),m,w) =y\right) \nonumber\\
&  \hspace{.5cm} \Cdot m(x) \Cdot \left[\SUMI{w^i_t=w}\right]\Bigg].
\end{align}
One can then write
\begin{align}
&\Prob{m_{t+1}(y)|m_t,\gamma_t} \substack{(a) \\ =} \Prob{\SUMI{x^i_{t+1}=y}}\\
&\substack{(b) \\ =} \Prob{\SUMI{f_t(x^i_t,\gamma(x^i_t),m_t,w^i_t)=y }}\\
&= \mathbb{P}\big(\sum_{x \in \mathcal{X}} \sum_{w \in \mathcal{W}} \SUMI{f_t(x_t,\gamma(x_t),m_t,w_t)=y } \\
& \hspace{5cm} \cdot \mathbb{1}(x^i_t=x,w^i_t=w) \big)\\
&\substack{(c) \\ =} \sum_{x \in \mathcal{X}} \sum_{w \in \mathcal{W}} \Big[\mathbb{1}\left( f_{t}(x,\gamma_t(x),m_t,w) =y\right)\\
& \hspace{2cm} \Cdot m_t(x) \Cdot \Prob{\SUMI{w^i_t=w}}\Big]\\
&\substack{(d) \\ =} \Prob{\bar f_t(m_t,\gamma_t,\mathbf{w}_t)(y)|m_t,\gamma_t},
\end{align}
where in the above equalities $(a)$ follows from \eqref{eq:mean-field}; $(b)$ follows from \eqref{eq:model1}; $(c)$ follows from the fact that $\mathbf{w}_t$ is i.i.d. and independent of $\mathbf x_t$ (and hence of $m_t$) and $\gamma_t$, and $(d)$ follows from \eqref{eq:def-bar-f}.  $\hfill \blacksquare$
\end{Proof}

\begin{Lemma}\label{lemma:bar f-hat f}
For every $m \in \mathcal{M}_n$ and $\gamma \in \mathcal{G}$, 
\begin{equation}
\Exp{\|\bar f_t(m,\gamma,\mathbf w) - \hat f_t(m,\gamma) \|_{\infty}} \leq \mathcal{O}(\frac{1}{\sqrt{n}}),
\end{equation}
where the expectation is taken with respect to $\mathbf{w} \in \mathcal{W}^n$.
\end{Lemma}
\begin{Proof}
For every $y \in \mathcal{X}$,  one has
\begin{align*}
&\Exp{\big|\bar f_t(m, \gamma,\mathbf w_t)(y)- \hat f_t(m,\gamma)(y)\big|}\\
&\qquad \substack{(a)\\=} \mathbb{E}\Big| \sum_{x \in \mathcal{X}} \sum_{w \in \mathcal{W}} m(x) \Cdot \mathbb{1} \left(f_t(x,\gamma(x),m,w)=y \right) \\
& \quad \qquad  \quad \Cdot \Big[ \SUMI{w^i_t=w} -\Prob{w} \Big]\Big| \\
&\qquad \substack{(b) \\ \leq } \sum_{x \in \mathcal{X}} \sum_{w \in \mathcal{W}} m(x) \Cdot \mathbb{1} \left(f_t(x,\gamma(x),m,w)=y \right) \\
& \quad \qquad  \quad \Cdot\mathbb{E}\Bigg|  \left[ \SUMI{w^i_t=w} -\Prob{w} \right]\Bigg| \substack{(c) \\ \leq } \mathcal{O}(\frac{1}{\sqrt{n}}),
\end{align*}
where  $(a)$ follows from \eqref{eq:relation-models}, \eqref{eq:mean-field-f}, and \eqref{eq:def-bar-f}; $(b)$ follows from the triangle inequality and monotonicity of the expectation operator, and $(c)$ follows from Lemma \ref{lemma:square-root-iid}, the fact that $\mathcal{X}$ and $\mathcal{W}$ are finite-valued spaces, and noting that $m(x) \leq 1$. Since the above result holds for every $y \in \mathcal{X}$,  it also holds in the infinity norm. $\hfill \blacksquare$
\end{Proof}
\B{Define the following dynamic program:
\begin{equation}\label{eq:value-function-Vz-T}
\hat V_{T+1}(z_{T+1}):=0, \quad  z_{T+1} \in \Delta(\mathcal{X}),
\end{equation}
and for  $z_t \in \Delta(\mathcal{X})$, $t=T,\ldots,1$,   
\begin{equation}\label{eq:value-function-Vz}
\hat V_t(z_t):=\min_{\gamma_t \in \mathcal{G}} \left(\hat c_t(z_t, \gamma_t)+\hat V_{t+1}\left( z_{t+1}\right) \right),
\end{equation}
where  $z_1(x)=\Prob{x_1=x}, x \in \mathcal{X}$, and  $z_t$ evolves deterministically as follows:
\begin{equation}\label{eq:z-dynamics}
 z_{t+1}= \hat{f}_t(z_t,\gamma_t).
\end{equation}}
\begin{Lemma}\label{lemma:dynamics- lipschitz-approximate-bound}
Let Assumption~\ref{assumption: finite-lipschitz} be satisfied. Then,  given $\gamma_t \in \mathcal{G}$, $m_t \in \mathcal{M}_n$, and $z_t \in \Delta(\mathcal{X}$),  the following inequality holds for $t \in \mathbb{N}_t$,
\begin{equation}
\Exp{\| m_{t+1} - z_{t+1} \|_{\infty}} \leq K^3_t \| m_t -z_t\|_{\infty}+\mathcal{O}(\frac{1}{\sqrt{n}}) .
\end{equation}
\end{Lemma}
\begin{Proof}
It is straightforward to show that
\begin{align}
& \Exp{\| m_{t+1} - z_{t+1} \|_{\infty}}
\substack{(a) \\ =} \Exp{\| \bar f_{t}(m_t,\gamma_t,\mathbf{w}_t) -  \hat f_t(z_t,\hspace{-.05cm}\gamma_t) \|_{\infty}}\\
& =  \Exp{\| \bar f_{t}(m_t,\gamma_t,\mathbf{w}_t) - \hat f_t(m_t,\gamma_t) \hspace{-.1cm}+ \hspace{-.1cm} \hat f_t(m_t,\gamma_t) \hspace{-.1cm}- \hspace{-.1cm} \hat f_t(z_t,\gamma_t) \|_{\infty}}\\ 
& \substack{(b) \\ \leq }  \mathbb{E}\| \hat f_t(m_t, \hspace{-.05cm}\gamma_t) \hspace{-.1cm}- \hspace{-.1cm} \hat f_t(z_t,\hspace{-.05cm} \gamma_t) \|_{\infty}\hspace{-.1cm}+\hspace{-.1cm}\mathbb{E}\| \bar f_{t}(m_t,\hspace{-.05cm}\gamma_t,\hspace{-.05cm} \mathbf{w}_t)\hspace{-.1cm}  - \hspace{-.1cm} \hat f_t(m_t,\hspace{-.05cm} \gamma_t)\|_{\infty} \hspace{-.1cm}\\
& \substack{(c) \\ \leq  }  K^3_t \Exp{\|m_t -z_t \|_{\infty}}+\mathcal{O}(\frac{1}{\sqrt{n}}), 
\end{align}
where $(a)$ follows from Lemma~\ref{lemma:bar f-mean-field} and \eqref{eq:z-dynamics}, $(b)$ follows from the triangle inequality, and $(c)$ follows from  Lemmas~\ref{lemma:Lipschitz} and~\ref{lemma:bar f-hat f}.~$\hfill \blacksquare$
\end{Proof}

Under mean-field sharing,  the  optimal solution for the system  described in Section~\ref{sec:problem-formulation} is identified by the following dynamic program~\cite[Theorem 1]{JalalCDC2014}:
\begin{equation} \label{eq:value-function-Vm-T}
V_{T+1}(m_{T+1}):=0, \quad  m_{T+1} \in \mathcal{M}_n,
\end{equation}
and  for  $m_t \in \mathcal{M}_n$,   $t=T,\ldots,1$, 
\begin{equation}\label{eq:value-function-Vm}
 V_t(m_t):=\min_{\gamma_t \in \mathcal{G}} \left(\hat c_t(m_t, \gamma_t)+ \Exp{V_{t+1}(m_{t+1})|m_t,\gamma_t} \right).
\end{equation}

\begin{Lemma}\label{lemma: value functions-difference}
Let Assumption~\ref{assumption: finite-lipschitz} hold. Then, given  $m_t \in \mathcal{M}_n$ and $z_t \in \Delta(\mathcal{X})$, at any time $t \in \mathbb{N}_T$, there exists a constant $K^5_t \in \mathbb{R}_{>0}$ such that  
\begin{equation}
\big|V_t(m_t) - \hat V_t(z_t) \big| \leq K^5_t \|m_t -z_t \|_{\infty} + \mathcal{O}(\frac{1}{\sqrt{n}}).
\end{equation}
\end{Lemma}
\begin{Proof}
The proof is presented in Appendix~\ref{sec:proof-lemma: value functions-difference}. $\hfill \blacksquare$
\end{Proof}

\begin{Theorem}\label{thm:finite-approximate-bound}
Let Assumption~\ref{assumption: finite-lipschitz} hold. Let also \B{$ \psi_t(z_t)$ be any argmin of the right-hand side of \eqref{eq:value-function-Vz} at time $t \in \mathbb{N}_T$. Define  fully decentralized strategy $\mathbf g:=\{g_t\}_{t=1}^{T}$  for  Problem~\ref{problem:finite}  such that}
\begin{equation}\label{eq:approximate-optimal-strategy}
g_t(x):=  \psi_t(z_t)(x), \quad x \in \mathcal{X}, t \in \mathbb{N}_T.
\end{equation}
Then, 
 \begin{equation}
\big| J(\mathbf{g}) - J^\ast \big| \leq \epsilon(n) \in \mathcal{O}(\frac{1}{\sqrt{n}}).
\end{equation}
\end{Theorem}
\begin{proof}
From the triangle inequality, one has
\begin{equation}\label{eq:error-performance-inequality}
| J(\mathbf g) - J^\ast| \leq  |J^\ast -  \hat V_1(z_1)|+| J(\mathbf g)  -  \hat V_1(z_1)|.
\end{equation}
It is desired now to show both terms on the right-hand side of \eqref{eq:error-performance-inequality} are $\mathcal{O}(\frac{1}{\sqrt{n}})$.

\textbf{Step 1:} In this step, we consider the first term of  \eqref{eq:error-performance-inequality}.  From \eqref{eq:value-function-Vz} and \eqref{eq:value-function-Vm} and noting that $J^\ast=\Exp{V_1(m_1)}$, we have
\begin{align}
&\big|J^\ast \hspace{-.1cm} -  \hspace{-.1cm}\hat V_1(z_1)\big| \hspace{-.1cm}= \hspace{-.1cm} \big| \Exp{V_1(m_1)} -\hat V_1(z_1)\big| \substack{a) \\ {=}}\big| \mathbb{E}[ V_1(m_1)-\hat V_1( z_1)]\big|\\
& \substack{(b) \\ {\leq}} \Exp{ |V_1(m_1) - \hat V_1(z_1) | }\substack{(c) \\ {\leq}} K^5_1\Exp{\|m_1 -z_1\|_{\infty}}+ \mathcal{O}(\frac{1}{\sqrt{n}})\\
&\substack{(d) \\ {\leq}} K^5_1 \mathcal{O}(\frac{1}{\sqrt{n}})+ \mathcal{O}(\frac{1}{\sqrt{n}}) = \mathcal{O}(\frac{1}{\sqrt{n}}),
\end{align}
where $(a)$ follows from the fact that $\hat V_1(z_1)$  is  deterministic; $(b)$ follows from  the monotonicity of the expectation operator,  which implies that   $\Exp{y} \leq \Exp{|y|}$   for every random variable $y$; $(c)$ follows from Lemma~\ref{lemma: value functions-difference}, and $(d)$ follows from  \eqref{eq:mean-field}, $z_1=\Prob{x_1}$, Lemma~\ref{eq:square-root-iid}, and the fact that the initial states are assumed to be i.i.d. random variables.

\B{\textbf{Step 2:} In this step, we consider the second term of  \eqref{eq:error-performance-inequality}.} Let $\hat m_t$ denote the empirical distribution of $n$ agents when they use strategy $\mathbf{g}$, given by Theorem~\ref{thm:finite-approximate-bound}, i.e., $\gamma_t=\psi_t(z_t)$. Therefore,
\begin{align}\label{eq:hat j- hat v- 1}
&\big| J(\mathbf g) -\hat V_1(z_1) \big| \substack{(a) \\ =}  \Big| \mathbb{E}\sum_{t=1}^T \hat c_t(\hat m_t,\psi_t(z_t)) - \sum_{t=1}^T \hat c_t(z_t,\psi_t(z_t))\Big| \nonumber\\
&\substack{(b) \\ \leq}  \sum_{t=1}^T K^4_t\mathbb{E} \| \hat m_t -z_t \|_{\infty},
\end{align}
where $(a)$ follows from \eqref{eq:total cost} and   \eqref{eq:value-function-Vz} (where $\min$ becomes equality under $\psi_t(z_t)$) and $(b)$ follows from Lemma~\ref{lemma:Lipschitz}. Note that $t=1$, $\hat m_1=m_1$, and for $t \in \mathbb{N}_{T-1}$, the evolution of the mean-field is given by Lemma~\ref{lemma:bar f-mean-field} as follows:
\begin{equation}\label{eq:hat m -dynamics}
\hat m_{t+1}=\bar f_t(\hat m_t, \psi_t(z_t), \mathbf{w}_t).
\end{equation}
Since both $\hat{m}_{t+1}$ (that evolves according to \eqref{eq:hat m -dynamics}) and $z_{t+1}$ (that evolves according to \eqref{eq:z-dynamics}) use identical strategy $\gamma_t=\psi_t(z_t)$, \B{we can use Lemma~\ref{lemma:dynamics- lipschitz-approximate-bound} to compute the expected difference}. In particular, given $\hat{m}_t$ and $z_t$, 
\begin{equation}\label{eq:hat j- hat v- 2}
\Exp{\| \hat m_{t+1} - z_{t+1} \|_{\infty}} \leq K^3_t \| \hat m_t -z_t\|_{\infty} +\mathcal{O}(\frac{1}{\sqrt{n}}) .
\end{equation}
Now, it results from  recursively  using  \eqref{eq:hat j- hat v- 2} in \eqref{eq:hat j- hat v- 1} and from  the monotonicity of the expectation operator that there exist  constants $K^6, K^7 \in \mathbb{R}_{>0}$ such that
\begin{multline}
| J(\mathbf g) -\hat V_1(z_1) | \leq K^6 \mathbb{E} \| \hat m_1 -z_1\|_{\infty} + K^7 \mathcal{O}(\frac{1}{\sqrt{n}}) \\
\substack{(c) \\=} K^6 \mathbb{E} \| m_1 -z_1\|_{\infty} + K^7 \mathcal{O}(\frac{1}{\sqrt{n}}) \substack{(d) \\ \leq} \quad  \mathcal{O}(\frac{1}{\sqrt{n}}),
\end{multline}
where $(c)$ follows from the fact that $\hat m_1=m_1$ and $(d)$ follows from Lemma~\ref{lemma:square-root-iid}, on noting that the initial states are i.i.d.  random variables. $\hfill \blacksquare$
\end{proof}

\begin{remark}
An important feature of  Theorem~\ref{thm:finite-approximate-bound} is that its proposed solution is independent of the number of agents $n$ because the functions $\hat f_t(\cdot)$ and $\hat c_t(\cdot)$,  given respectively by \eqref{eq:mean-field-f} and \eqref{eq:mean-field-cost}, are independent of $n$.
\end{remark} 
\begin{remark}
Note that \eqref{eq:value-function-Vz-T}, \eqref{eq:value-function-Vz} and \eqref{eq:z-dynamics} do not depend on the information of agents; hence, they may be solved off-line.  More precisely, every agent can independently compute   $\{\psi_t(z_t)\}_{t=1}^{T}$ in a distributed manner with no communication required. In the case of multiple solutions, agents can make sure they all  compute the same solution while using $\argmin$ by agreeing upon a deterministic rule to break a tie.  Thus, strategy \eqref{eq:approximate-optimal-strategy} can be implemented  based on a completely decentralized information structure \eqref{eq:no-sharing}.
\end{remark}
 Initially, every agent locally computes  the $\argmin$ of \eqref{eq:value-function-Vz} for all $t \in \mathbb{N}_T$, i.e., $\{\psi_t(z_t)\}_{t=1}^{T}$. Then,  when the system is operating, agent $i$ makes a decision based on the local state $x^i_t$, i.e., 
\[u^i_t=\psi_t(z_t)(x^i_t).\]
According to Theorem~\ref{thm:finite-approximate-bound}, the  decision of each agent $i \in \mathbb{N}_n$ at time $t \in \mathbb{N}_T$ is determined by three factors: (a)  strategy $\psi_t$ that depends on the model of the system; (b)  variable  $z_t$ that is common knowledge among all agents and evolves to  $z_{t+1}$ according to \eqref{eq:z-dynamics},  and (c)  local state $x^i_t$   that is only known to agent $i$.

\begin{Proposition}\label{cor:continuous-lip}
Let  Assumption~\ref{assumption: finite-lipschitz} hold. Then,  $\hat V_t(z)$ is a Lipschitz function, i.e., there exists a constant $K_t \in \mathbb{R}_{ >0}$, $t \in \mathbb{N}_T$,  such that for every $z_1,z_2 \in \mathcal{I}(\mathcal{X})$,
\begin{equation}
\big|  \hat V_t (z_1) - \hat V_t(z_2) \big| \leq K_t \| z_1 - z_2\|_{\infty}.
\end{equation}
\end{Proposition}
\begin{Proof}
The proof is  presented in Appendix~\ref{sec:proof-cor:continuous-lip}. $\hfill \blacksquare$
\end{Proof}

The recursion introduced in \eqref{eq:value-function-Vz-T}--\eqref{eq:z-dynamics} is computationally intractable, in general, because $\Delta(\mathcal{X})$ is an uncountable set. 
%
 According to Lemma~\ref{lemma:Lipschitz} and Proposition~\ref{cor:continuous-lip},  $\hat{f}(z)$ and  $\hat{V}(z)$ are Lipschitz continuous, respectively.  Therefore, one could quantize the  infinite-set optimization of Theorem~\ref{thm:finite-approximate-bound} into a finite-set one such that the quantization  error  is upper-bounded  by some Lipschitz function. Using this idea,  it is shown  in Corollary~\ref{cor:quatize} that  the rate of convergence in Theorem~\ref{thm:finite-approximate-bound} is preserved under a uniform  quantization.
\begin{Corollary}\label{cor:quatize}
Let Assumption~\ref{assumption: finite-lipschitz} hold and  function $Q: \Delta(\mathcal{X}) \rightarrow \mathcal{Q}_n$  map every point $z \in \Delta(\mathcal{X})$ to its nearest point $\hat z$ in $\mathcal{Q}_n$, i.e.,
\begin{equation}
Q(z) \in \argmin_{\hat z \in \mathcal{Q}_n} \|z -\hat z \|_{\infty}.
\end{equation}
Define
\begin{equation}
 \hat V_{T+1}(\hat z_{T+1}):=0, \quad  \hat z_{T+1} \in \mathcal{Q}_n,
\end{equation}
and for  $t=T,\ldots,1$, $\hat z_t \in \mathcal{Q}_n$,
\begin{equation}\label{eq:value-function-Vz-cor}
\hat V_t(\hat z_t):=\min_{\gamma_t \in \mathcal{G}} \left(\hat c_t(\hat z_t, \gamma_t)+ \hat V_{t+1}\left(\hat z_{t+1}\right) \right),
\end{equation}
where  $\hat z_1=Q(\Prob{x_1})$, $\Prob{x_1} \in \Delta(\mathcal{X}),$ and  $\hat z_t$ evolves deterministically as follows
\begin{equation}
 \hat z_{t+1}= Q(\hat{f}_t(\hat z_t,\gamma_t)).
\end{equation}
 Let $ \psi_t(\hat z_t)$ be any argmin of the right-hand side of \eqref{eq:value-function-Vz-cor} and define  $\mathbf g:=\{g_t\}_{t=1}^{T}$,   where
\begin{equation}
g_t(x):=  \psi_t(\hat z_t)(x), \quad x \in \mathcal{X}, t \in \mathbb{N}_T.
\end{equation}
Then, strategy $\mathbf{g}$ is a solution to Problem~\ref{problem:finite} such that
 \begin{equation}
\big| J(\mathbf{g}) - J^\ast \big| \leq \epsilon(n) \in \mathcal{O}(\frac{1}{\sqrt{n}}).
\end{equation}
\end{Corollary}
\begin{Proof}
The error associated with quantizing $z \in \Delta(\mathcal{X})$ into $\hat z \in \mathcal{Q}_n$ is  bounded by  $\sup_{z \in \Delta(\mathcal{X})}\|z -Q(z) \|_{\infty} \leq \frac{1}{n}$.
According to  Lemma~\ref{lemma:Lipschitz} and Proposition~\ref{cor:continuous-lip}, functions $\hat f_t$ and $\hat V_t$ are Lipschitz in $z$. Therefore, this quantization  error will be upper-bounded by $\mathcal{O}(\frac{1}{n})$ over a fixed finite horizon $T$.  The proof is complete on nothing that  $\mathcal{O}(\frac{1}{n})$  is dominated by  $\mathcal{O}(\frac{1}{\sqrt{n}})$,  when $n$ is large.  $\hfill \blacksquare$
\end{Proof}

\section{Conclusions}\label{sec:conclusion}
In this paper,  team-optimal control of a large number of homogeneous agents modeled as mean-field coupled Markov chains is considered. Every agent observes only its own local state, i.e., the strategy has a fully decentralized information structure.  A sub-optimal strategy, independent of the number of agents~$n$,  is proposed whose performance converges to that of the optimal  mean-field sharing strategy at the rate  $1 / \sqrt{n}$.  To establish this result,  it is assumed that the transition probabilities and costs  are Lipschitz continuous in  the mean-field (i.e., no assumption is imposed on the strategy). To find the sub-optimal strategy, it is required to solve an infinite-set optimization problem, in general. To address this concern,  a novel idea is proposed  to quantize the  infinite-set optimization problem into a finite-set one.  In particular, it is shown that under uniform quantization with the step-size of $1/n$,  the convergence rate $1 / \sqrt{n}$ is preserved.  
\bibliography{CDC2017}
\bibliographystyle{IEEEtran}

\section*{Appendix}

\subsection{Proof of Lemma~\ref{lemma: value functions-difference}}\label{sec:proof-lemma: value functions-difference}
We use backward induction to prove this lemma. At $t=T$, \eqref{eq:value-function-Vm-T} and \eqref{eq:value-function-Vm} yield
\begin{align}
V_T(m_T)&=\min_{\gamma_t \in \mathcal{G}} \hat{c}_T(m_T,\gamma_T)\\
& \substack{(a) \\ \leq } \min_{\gamma_t \in \mathcal{G}} \left(\Big|  \hat{c}_T(m_T,\gamma_T) - \hat c_t(z_T,\gamma_T)\Big| + \hat c_t(z_T,\gamma_T)\right)\\
& \substack{(b) \\ \leq } \min_{\gamma_t \in \mathcal{G}} \left(  K^4_T \| m_T -z_T\|_{\infty}+ \hat c_t(z_T,\gamma_T)\right)\\
& =   K^4_T \| m_T -z_T\|_{\infty}+ \min_{\gamma_t \in \mathcal{G}}  \hat c_t(z_T,\gamma_T)\\
&\substack{(c) \\ \leq } K^4_T \| m_T -z_T\|_{\infty}+ \hat V_T(z_T),
\end{align}
where $(a)$ follows from the triangle inequality,  per-step costs being non-negative (by definition), and the monotonicity of the minimum operator, $(b)$ follows from Lemma~\ref{lemma:Lipschitz} and the monotonicity of minimum operator, and $(c)$ follows from~\eqref{eq:value-function-Vz}. Therefore, $K^5_T:=K^4_T$, i.e.,
\begin{equation}
\Big| V_T(m_T) -\hat{V}_T(z_T) \Big| \leq K^5_T \| m_T -z_T\|_{\infty}.
\end{equation}
Suppose now that the inequality holds at time $t+1$, i.e., 
\begin{equation}\label{eq:indunction-t+1-approximate}
\Big| V_{t+1}(m_{t+1}) -\hat{V}_{t+1}(z_{t+1}) \Big| \leq K^5_{t+1} \| m_{t+1} -z_{t+1}\|_{\infty} + \mathcal{O}(\frac{1}{\sqrt{n}}).
\end{equation}
It is desired to  prove the result for time $t$.  It is deduced from \eqref{eq:value-function-Vm} that
\begin{align*}
V_t(m_t)&=\min_{\gamma_t \in \mathcal{G}} \left( \hat c_t(m_t,\gamma_t)+ \Exp{V_{t+1}(m_{t+1})|m_t,\gamma_t} \right)\\
&=  \min_{\gamma_t \in \mathcal{G}} \Big(\hat c_t(m_t,\gamma_t) - \hat c_t(z_t,\gamma_t)  \\
&\quad + \Exp{V_{t+1}(m_{t+1})|m_t,\gamma_t}-  \hat V_{t+1}(z_{t+1})\\
&\quad+\hat c_t(z_t,\gamma_t)+ \hat V_{t+1}\left(z_{t+1}\right)\Big)\\
&\substack{(d) \\ \leq } \min_{\gamma_t \in \mathcal{G}} \Big( \Big| \hat c_t(m_t,\gamma_t) - \hat c_t(z_t,\gamma_t) \Big| \\
&\quad + \Exp{\Big| V_{t+1}\left(m_{t+1}\right) -  \hat V_{t+1}\left(z_{t+1}\right) \Big| | m_t,z_t,\gamma_t}\\
&\quad + \hat c_t(z_t,\gamma_t)+ \hat V_{t+1}\left(z_{t+1}\right)\Big)\\
&\substack{(e) \\ \leq } \min_{\gamma_t \in \mathcal{G}} \Big(K^4_t \| m_t-z_t\|_{\infty}\\
&\quad  +K^5_{t+1}\Exp{ \| m_{t+1}-z_{t+1}\|_{\infty}|m_t,z_t,\gamma_t} + \mathcal{O}(\frac{1}{\sqrt{n}})\\
&\quad+\hat c_t(z_t,\gamma_t)+ \hat V_{t+1}\left(z_{t+1}\right)\Big)\\
&\substack{(f) \\ \leq } (K^4_t\hspace{-.1cm}+\hspace{-.1cm}K^5_{t+1}K^3_t) \| m_t-z_t\|_{\infty} + (1+K^5_{t+1})\mathcal{O}(\frac{1}{\sqrt{n}})\\
&\quad + \min_{\gamma_t \in \mathcal{G}}\Big(\hspace{-.05cm}\hat c_t(z_t,\gamma_t) \hspace{-.1cm}+\hspace{-.1cm} \hat V_{t+1}\left(z_{t+1}\right)\hspace{-.1cm}\Big)\\
&\substack{(g) \\= } (K^4_t\hspace{-.1cm}+\hspace{-.1cm}K^5_{t+1}K^3_t) \| m_t-z_t\|_{\infty} +\mathcal{O}(\frac{1}{\sqrt{n}}) + \hat V_t(z_t),
\end{align*}
where $(d)$  follows from the triangle inequality,  per-step costs being non-negative (by definition), and the monotonicity of minimum and expectation operators; $(e)$ follows from Lemma~\ref{lemma:Lipschitz} and \eqref{eq:indunction-t+1-approximate};  $(f)$  follows from Lemma~\ref{lemma:dynamics- lipschitz-approximate-bound}, and $(g)$  follows from \eqref{eq:value-function-Vz}. Therefore, there exists a constant $K^5_t:= K^4_t+ K^5_{t+1}K^3_t$ such that
\begin{equation}\label{eq:indunction-t-approximate}
\Big| V_{t}(m_{t}) -\hat{V}_{t}(z_{t}) \Big| \leq K_{t} \| m_{t} -z_{t}\|_{\infty} + \mathcal{O}(\frac{1}{\sqrt{n}}).
\end{equation}

\subsection{Proof of Proposition~\ref{cor:continuous-lip}}\label{sec:proof-cor:continuous-lip}
We use backward induction. At $t=T$, 
\begin{align*}
\hat{V}_T(z_1)&= \min_{\gamma_T} \hat{c}_T(z_1,\gamma_T)\\
&=\min_{\gamma_T} \left( \hat{c}_T(z_1,\gamma_T) -\hat{c}_T(z_2,\gamma_T)+\hat{c}_T(z_2,\gamma_T) \right)\\
& \substack{(a) \\ \leq } \min_{\gamma_T} \left( \big|\hat{c}_T(z_1,\gamma_T) -\hat{c}_T(z_2,\gamma_T) \big|+\hat{c}_T(z_2,\gamma_T) \right)\\
& \substack{(b) \\ \leq } \min_{\gamma_T} \left( K^4_T \| z_1-z_2\|_{\infty}+\hat{c}_T(z_2,\gamma_T) \right)\\
&= K^4_T \| z_1-z_2\|_{\infty}+ \min_{\gamma_T} \hat{c}_T(z_2,\gamma_T) \\
&\substack{(c) \\ =} K^4_T \| z_1-z_2\|_{\infty}+ \hat V_T(z_2),
\end{align*}
where $(a)$ follows from the triangle inequality,  the fact that $\hat c_T(\Cdot) \in \mathbb{R}_{\geq 0}$, and  the monotonicity of the minimum operator; $(b)$ follows from Lemma~\ref{lemma:Lipschitz} and the monotonicity of the minimum operator, and $(c)$ follows from~\eqref{eq:value-function-Vz}.
Suppose now that the inequality holds at time $t+1$, i.e.,  
\begin{equation}\label{eq:t+1-step-induction-hat-v}
\big|\hat{V}_{t+1}(z_1) - \hat{V}_{t+1}(z_2)\big| \leq  K_{t+1} \| z_1-z_2\|_{\infty}.
\end{equation}
It is desired to  prove the the result for time $t$. One can write
\begin{align*}
\hat V_t(z_1)&= \min_{\gamma_t} \big(\hat c_t(z_1,\gamma_t)+\hat{V}_{t+1}(\hat{f}_t(z_1,\gamma_t))\big)\\
&=\min_{\gamma_t}\Big(  \hat c_t(z_1,\gamma_t) - \hat c_t(z_2,\gamma_t)+ \hat c_t(z_2,\gamma_t)\\
& +  \hat{V}_{t+1}(\hat{f}_t(z_1,\gamma_t)) \hspace{-.1cm}- \hspace{-.1cm}\hat{V}_{t+1}(\hat{f}_t(z_2,\gamma_t)) \hspace{-.1cm}+ \hspace{-.1cm}\hat{V}_{t+1}(\hat{f}_t(z_2,\gamma_t)) \Big)\\
&\substack{(d) \\ \leq} \min_{\gamma_t}\Big( \big| \hat c_t(z_1,\gamma_t) - \hat c_t(z_2,\gamma_t)\big| \\
& \quad + \big| \hat{V}_{t+1}(\hat{f}_t(z_1,\gamma_t))-\hat{V}_{t+1}(\hat{f}_t(z_2,\gamma_t)) \big|\\
& \quad +\hat c_t(z_2,\gamma_t)+\hat{V}_{t+1}(\hat{f}_t(z_2,\gamma_t)) \Big)\\
&\substack{(e) \\ \leq}(K^4_t+K_{t+1}K^3_t)\|z_1 -z_2 \|_{\infty} \\
&\quad + \min_{\gamma_t} \big(\hat c_t(z_2,\gamma_t)+\hat{V}_{t+1}(\hat{f}_t(z_2,\gamma_t))\big)\\
&\substack{(f) \\ =}(K^4_t+K_{t+1}K^3_t)\|z_1 -z_2 \|_{\infty} +   \hat V_t(z_2),
\end{align*}
where $(d)$ follows from the triangle inequality,  the fact that $\hat c_t (\Cdot), \hat V_{t+1}(\Cdot) \in \mathbb{R}_{ \geq 0}$, and the monotonicity of the minimum operator; $(e)$ follows from Lemma~\ref{lemma:Lipschitz}, \eqref{eq:t+1-step-induction-hat-v}, and  the monotonicity of the minimum operator, and $(f)$ follows from~\eqref{eq:value-function-Vz}. $\hfill \blacksquare$

\end{document}